\documentclass[12pt]{article}

\newtheorem{thm}{Theorem}[section]

\newtheorem{proposition}[thm]{Proposition}

\newenvironment{proof}{\par\noindent Proof:\ }{${\tt [}\kern-0.2mm{\tt ]}$}

\newcommand{\dx}{\partial_x}
\newcommand{\dy}{\partial_y}
\newcommand{\dz}{\partial_z}
\newcommand{\Ext}{{\rm Ext}}
\newcommand{\gr}{{\rm gr}}
\newcommand{\lra}{\longrightarrow}

\newcommand{\Ann}{Ann_{\mathcal D}(1/f)}
\newcommand{\DD}{{\mathcal D}}
\newcommand{\cO}{{\mathcal O}}
\newcommand{\vv}{\vspace{1cm}}
\newcommand{\Der}{{{\rm Der(log} \,f)}}

\date{June 30, 2000}
\title{On some ${\cal D}$-modules in dimension
2}
\author{Francisco Jes\'us Castro Jim\'enez, Jos\'e
Mar\'{\i}a Ucha Enr\'{\i}quez\thanks{Both authors partially supported by
DGESIC-PB97-0723 and by A.I. HF1998-0105}. \\ Facultad de
Matem\'{a}ticas. Apdo. 1160. E-41080 SEVILLA (SPAIN)\\ castro@cica.es;
ucha@algebra.us.es}
\begin{document}

\maketitle

{\bf Abstract.-} We prove a duality formula for two ${\cal
D}$-modules arising from logarithmic derivations w.r.t. a plane
curve. As an application we give a differential proof of a
logarithmic comparison theorem in \cite{4T}.

{\bf Math.
Classification:} 32C38, 13N10, 14F40, 13P10.

\section{Introduction}

Let ${\mathcal O}={\bf C}\{x,y\}$ be the ring of convergent power
series in two variables and $\mathcal D$ the ring of linear
differential operators with coefficients in $\mathcal O$. For each
reduced  power series $f\in \mathcal O$, with $f(0,0)=0$,  we will
denote by $I^{log}$ the left ideal of $\mathcal D$ generated by
the logarithmic derivations (see \cite{Saito}) with respect to
$f$. We denote by ${\rm Der}_{\bf C}(\cO)$ the Lie algebra of
${\bf C}$-derivations on $\cO$. Recall that a derivation
$\delta\in {\rm Der}_{\bf C}({\mathcal O})$ is logarithmic if
there exists $a\in {\mathcal O}$ such that $\delta(f)=af$. We
denote by $\widetilde{I}^{log}$ the left ideal of ${\mathcal D}$
generated by the operators of the form $\delta +a$ where
$\delta(f)=af$.

We  first prove that the ${\mathcal D}$-modules $M^{log}={\mathcal
D}/I^{log}$ and $\widetilde{M}^{log}={\mathcal
D}/\widetilde{I}^{log}$ are dual each to the other and then that
both ${\mathcal D}$-modules are regular holonomic \ref{dual}.

Let ${\mathcal O}[1/f]$ be the ${\mathcal D}$-module of (the germs
of) the meromorphic functions in two variables with poles along
$f$. There exists a natural surjective morphism $\psi :
\widetilde{M}^{log} \rightarrow {\mathcal O}[1/f]$. Using
\cite{4T} we prove that $\Ext_{{\mathcal
D}}^2(\widetilde{M}^{log},{\mathcal O})=0$ if and only if $f$ is
quasi-homogeneous and then we obtain that the morphism $\psi$ is
an isomorphism if and only if $f$ is quasi-homogeneous (see
\ref{nocero}). As a consequence we give a new ``differential"
proof of the logarithmic comparison theorem of \cite{4T}.

These results are susceptible to be generalized to the case of
higher dimensions but no general results are known up to now. See
\cite{Castro-Ucha-dualite} for a proof of the duality formula in
higher dimension. Nevertheless we give a complete example showing
that some results of the present work are true in dimension 3.

We wish to thank Prof. L. Narv\'aez for giving us useful
suggestions.

\section{The module $\widetilde{M}^{log}$ in the general case.}

Let us consider any reduced $f \in {\mathcal O} = {\bf C} \{ x, y
\}$ with a singular point at the origin. It is possible to obtain,
from the logarithmic derivations, an ideal inside $Ann_{\DD}
(1/f)$: if $\delta(f) = af$ then $\delta + a \in Ann_{\DD}(1/f)$.
This fact suggested us a general way to present the annihilating
ideal  of $1/f$ for a constructive proof of the equality
$\Ext^2_{\DD}({\mathcal O}[1/f], {\mathcal O}) = 0$ for any
``polynomial" curve (see \cite{Ucha}).

We have $\widetilde{I}^{\log} \subset \Ann$, where
$\widetilde{I}^{\log}$ is the left ideal in ${\mathcal D}$
generated by the operators $\delta+a$ for $\delta\in Der_{\bf
C}(\cO)$ and $\delta(f)=af$. Then we have a surjective morphism
$\psi: \widetilde{M}^{log}=\frac{{\mathcal
D}}{\widetilde{I}^{\log}} \longrightarrow \DD / \DD Ann_\DD (1/f)
\simeq {\mathcal O}[1/f]$ (for the last isomorphism we use that
the Bernstein polynomial of $f$ has no integer roots smaller than
-1 (see \cite{Va})). It is well known that around each smooth
point of $f=0$ the morphism $\psi$ is in fact an isomorphism. So,
the kernel $K$ of $\psi$ is a ${\mathcal D}$-module concentrated
at the origin. Then $K$  is a direct sum of ``couches-multiples"
modules \cite{K}, and this type of modules are regular holonomic
\cite{MK2}. In particular $\widetilde{M}^{\log}$ is regular
holonomic because ${\mathcal O}[\frac{1}{f}]$ and $K$ are.

\noindent We will denote by $\Der$ the Lie algebra of logarithmic
derivations with respect to $f$. By \cite{Saito} $\Der$ is a free
$\cO$-module of rank two. Let $\{ \delta_1, \delta_2 \}$ be a
basis of Der(log $f$), $$ \left\{
\begin{array}{ccl}
\delta_1 & = & b_1 \dx + c_1 \dy, \\ \delta_2 & = & b_2 \dx + c_2
\dy
\end{array}
\right. . $$ \noindent We can suppose that $$ \left|
\begin{array}{cc}
 b_1 & c_1 \\
 b_2 & c_2
\end{array}
\right|
=
f $$ \noindent We will take into account the following results for
any (reduced) curve $f$:
\begin{itemize}
\item Every basis $\delta_1, \delta_2$ of Der(log $f$)
verifies that $$ \langle \sigma(\delta_1), \sigma(\delta_2)
\rangle = \gr^F (I^{log}) = \gr^F (\widetilde{I}^{log}), $$
because $\{ \sigma(\delta_1) , \sigma(\delta_2) \}$ is a regular
sequence (see \cite{Cald} and (\cite{Cald3}, Corollary 4.2.2)).
Here $\sigma(\cdot)$ denotes the principal symbol of the
corresponding operator and $\gr^F (I^{log})$ is the graded ideal
associated to the order filtration on ${\mathcal D}$. Therefore,
$$ CCh (\widetilde{M}^{log}) = CCh (M^{log}), $$ where $CCh(\ )$
represents the {\em characteristic cycle} of the $\DD$-module
(see, for example, \cite{G-M}).  Of course both $M^{log}$ and
$\widetilde{M}^{log}$ define coherent $\mathcal D$-modules in some
neighborhood of the origin and then we can properly speak of
characteristic varieties and characteristic cycles. Since
$\widetilde{M}^{log}$ is holonomic then $M^{log}$ is holonomic.
\item For any curve,
$$ Sol (M^{log}) \stackrel{q.i.}\simeq \Omega^{\bullet} (log f)
\stackrel{\varphi}\lra \Omega^{\bullet}[1/f] \simeq DR ({\mathcal
O}[1/f]), $$ where $Sol(\ )$ and $DR(\ )$ are the {\em solutions
complex} and the {\em De Rham complex} (see, for example,
\cite{MK2}) and where $\Omega^\bullet(\log f)$ (resp.
$\Omega^\bullet([1/f])$) is the complex of logarithmic
differential forms (resp. meromorphic differential forms). The
first quasi-isomorphism appears in \cite{Cald3} and $\varphi$ is
the natural morphism.
\end{itemize}

\begin{proposition}{\label{alfa}}
Let $f$ be a (reduced) curve and let $\{ \delta_1, \delta_2 \}$ be
a basis of {\rm Der (log $f$)} with $[ \delta_1, \delta_2 ] =
\alpha_1 \delta_1 + \alpha_2 \delta_2$ and  $\delta_i (f) = a_i
f,\ i=1,2$. Then $$ {\mathcal D}\{\delta_2^t + \alpha_1 ,
\delta_1^t - \alpha_2 \} = {\mathcal D} \{\delta_1 + a_1 ,
\delta_2 + a_2 \} $$ \noindent where $\delta_i^t$ is the {\em
transposed} of $\delta_i$.
\end{proposition}
\begin{proof}
First we find an expression of the $\alpha_i$ from the $a_j, b_k,
c_l$: $$ [\delta_1, \delta_2] = \alpha_1 (b_1 \dx + c_1 \dy) +
\alpha_2 (b_2 \dx + c_2 \dy) = $$ $$ = (\alpha_1 b_1 + \alpha_2
b_2) \dx + (\alpha_1 c_1 + \alpha_2 c_2) \dy = $$ $$ = b_1 \dx
(b_2) \dx - b_2 \dx (b_1) \dx + b_1 \dx (c_2) \dy - c_2 \dy (b_1)
\dx + $$ $$ + c_1 \dy (b_2) \dx - b_2 \dx (c_1) \dy + c_1 \dy
(c_2) \dy - c_2 \dy (c_1) \dy = $$ $$ = (c_1 \dy (b_2) - b_2 \dx
(b_1) - c_2 \dy (b_1) + c_1 \dy (b_2)) \dx + $$ $$ + (b_1 \dx
(c_2) - b_2 \dx (c_1) - c_2 \dy (c_1) + c_1 \dy (c_2)) \dy. $$
\noindent Besides, $$ - \delta_1^t + \alpha_2 =  \dx b_1 + \dy c_1
+ \alpha_2=  \delta_1 + \alpha_2 + \dx (b_1) + \dy (c_1). $$
\noindent To prove that $\alpha_2 + \dx (b_1) + \dy (c_1) = a_1$,
we will establish that $$ \alpha_2 f = a_1 f - \dx (b_1) f - \dy
(c_1) f $$ \noindent We have $$ a_1 f - \dx (b_1) f - \dy (c_1) f
= $$ $$ = (b_1 \dx + c_1 \dy) - \dx (b_1) - \dy (c_1)) (b_1 c_2 -
b_2 c_1) = $$ $$ = b_1 (b_1 \dx (c_2) - c_1 \dx (b_2) - b_2 \dx
(c_1)) + $$ $$ + c_1 (c_2 \dy (b_1) + b_1 \dy (c_2) - c_1 \dy
(b_2)) + $$ $$ + b_2 c_1 \dx (b_1) - b_1 c_2 \dy (c_1). $$
Therefore $$ (\alpha_1, \alpha_2) \left(
\begin{array}{cc}
b_1 & c_1 \\ b_2 & c_2
\end{array}
\right) = (\gamma_1, \gamma_2), $$ \noindent where $$
\begin{array}{ccl}
\gamma_1 & = & c_1 \dy (b_2) - b_2 \dx (b_1) - c_2 \dy (b_1) + c_1
\dy (b_2), \\ \gamma_2 & = & b_1 \dx (c_2) - b_2 \dx (c_1) - c_2
\dy (c_1) + c_1 \dy (c_2)).
\end{array}
$$ \noindent Multiplying by the transposed adjoint matrix and by
$f$ we obtain $$ (\alpha_1 f, \alpha_2 f) = (\gamma_1, \gamma_2)
\left(
\begin{array}{cc}
c_2 & - c_1 \\ - b_2 & b_1
\end{array}
\right) . $$ \noindent and hence the equality follows. In a
similar way $\delta_2^t + \alpha_1 = - \delta_2 - a_2$. Then both
ideals are equal.
\end{proof}

\vv

Prof. Narv\'aez pointed us to consider, instead of the Lie algebra
Der(log$f$), the Lie algebra $$ L = \{ \delta + a |\ \delta(f) =
af \}, $$ \noindent and try to construct of a free resolution (of
``Spencer type") of $\widetilde{M}^{log}$ (\cite{Cald},
\cite{Cald3}). In fact, we have

\begin{proposition}{\label{resol}}
A free resolution of $\widetilde{M}^{log}$ is $$ 0 \lra \DD
\stackrel{\phi_2} \lra \DD^2 \stackrel{\phi_1}\lra \DD \lra
\widetilde{M}^{log} \lra 0, $$ \noindent where $\phi_2$ is defined
by the matrix $$ (- \delta_2 - a_2 - \alpha_1,  \delta_1 + a_1 -
\alpha_2), $$ \noindent and $\phi_1$ by $\left( \begin{array}{c}
\delta_1 + a_1 \\ \delta_2 + a_2
\end{array}\right)$.
\end{proposition}
\begin{proof}
To check the exactness of the resolution above, it is enough to
consider a discrete filtration on that complex and to verify the
exactness of the resulting resolution (see \cite{B}, chapter 2,
lemma 3.13). The same argument is used in \cite{Cald},
\cite{Cald3}(proposition 4.1.3) to prove that the complex $\DD
\otimes_{{\mathcal V}_0^f (\DD)} Sp^{\bullet} (log f)$ is a free
resolution of $M^{log}$ (as a left $\DD$-module)\footnote{Here
$Sp^{\bullet} (log f)$ is the Logarithmic Spencer complex and
${\mathcal V}_0^f (\DD)$ is the ring of degree zero differential
operators w.r.t. $\cal V$-filtration relative to $f$. See
\cite{Cald}, \cite{Cald3}, section 1.2.}. But, for $n=2$, the
exact graded complex in the proof of \cite{Cald3} is precisely $$
0 \lra \gr^F(\DD) \stackrel{M_1}\lra \gr^F(\DD)^2[-2]
\stackrel{M_2}\lra \gr^F (\DD) [-1] \lra \gr^F (M^{log}) \lra 0,
$$ \noindent where the matrices are $$ M_1 = (-\sigma^F(\delta_2),
\sigma^F(\delta_1)),  \ \ M_2 = \left(
\begin{array}{c} \sigma^F(\delta_1) \\ \sigma^F(\delta_2)
\end{array} \right). $$ \noindent And the last complex is the
result of applying the same graduation to the resolution of
$\widetilde{M}^{log}$ too, because $$ \sigma^F (\delta_i) =
\sigma^F (\delta_i + a_i). $$
\end{proof}

\begin{proposition}{\label{dual}}
Given $f \in {\bf C} \{x,y\}$ , $\widetilde{M}^{log} \simeq
(M^{log})^\star$ where $({ })^\star$ is the dual in the sense of
$\cal D$-modules. In particular $\widetilde{M}^{log}$ and
$M^{log}$ are regular $\DD$-modules.
\end{proposition}
\begin{proof}
We take the free resolution of $M^{log}$ (see \cite{Cald},
(\cite[Th. 3.1.2]{Cald3}) $$ 0 \lra \DD \stackrel{\psi_2}\lra
\DD^2 \stackrel{\psi_1}\lra \DD \lra M^{log} \lra 0, $$ \noindent
where $\{ \delta_1, \delta_2 \}$ is a basis of the ${\mathcal
O}$-module Der(log$f$), where $$ [ \delta_1, \delta_2 ] = \alpha_1
\delta_1 + \alpha_2 \delta_2, $$ $$ \psi_1 =
{\left(\begin{array}{c} \delta_1
\\ \delta_2
\end{array}\right)}
$$ \noindent and, on the other hand, $\psi_2$ is the syzygy matrix
$$ \psi_2 = (- \delta_2 - \alpha_1, \delta_1 - \alpha_2). $$
\noindent Applying the $Hom_{\DD}( -, \DD)$ functor to calculate
the dual module, we obtain the sequence $$ 0 \lra \DD
\stackrel{\psi_1^*} \lra \DD^2 \stackrel{\psi_2^*}\lra \DD \lra 0,
$$ \noindent where $\psi_2^*$ is the right product by
${\left(\begin{array}{c} - \delta_2 - \alpha_1 \\ \delta_1 -
\alpha_2
\end{array}\right)}$.
Hence, $(M^{log})^\star$ is the left $\DD$-module associated to
the right $\DD$-module $\DD/(\delta_2 + \alpha_1, \delta_1 -
\alpha_2) \DD$, that is to say, $$ (M^{log})^\star \simeq
\DD/\DD(\delta_2^t+ \alpha_1, \delta_1^t - \alpha_2). $$ \noindent
Using the proposition \ref{alfa}, we deduce that $(M^{log})^\star
\simeq \widetilde{M}^{log}$. The regularity of $M^{log}$ follows
from the regularity of $\widetilde{M}^{log}$ (c.f. \cite{MK2}).
\end{proof}

\begin{proposition}{\label{ext_no_cero}}
If $f$ is a non quasi homogeneous (reduced) curve, then $$
\Ext^2_{\DD} (\widetilde{M}^{log}, {\mathcal O}) \neq 0. $$
\end{proposition}
\begin{proof}
The proof of this result contains, as an essential ingredient, a
re-reading of the demonstration of Theorem 3.7 of \cite{4T}. As a
matter of fact, we include some tricks of this demonstration.

By proposition \ref{resol}, a free resolution of
$\widetilde{M}^{log}$ is $$ 0 \lra \DD \stackrel{\phi_2} \lra
\DD^2 \stackrel{\phi_1}\lra \DD \lra \widetilde{M}^{log} \lra 0,
$$ \noindent where $\phi_2$ is the matrix $$ (- \delta_2 - a_2 -
\alpha_1,  \delta_1 + a_1 - \alpha_2). $$ \noindent Hence,
$\Ext^2_{\DD}(\widetilde{M}^{log}, {\mathcal O}) \simeq {\mathcal
O}/Img \phi_2^*$. To guarantee that this vector space has
dimension greater than zero, it is enough to show that a pair of
functions $h_1, h_2 \in {\mathcal O}$ such that $$ (-\delta_2 -
a_2 - \alpha_1, \delta_1 + a_1 - \alpha_2) \left(
\begin{array}{c}
h_1 \\ h_2
\end{array}
\right)
=
1, $$ \noindent does not exist, that is to say, that $1 \notin Img
\phi_2^*$.

Let us take $\delta_1 = b_1 \dx + c_1 \dy$. As  $a_1 - \alpha_2 =
\dx (b_1) + \dy (c_1)$, (proposition \ref{alfa}) we will prove
that, or $b_1$ and $c_1$ have no lineal parts, or that after
derivation those lineal parts become 0.

Of course $f$ has no quadratic part: in that case, because of the
classification of the singularities in two variables, $f$ would be
equivalent to a quasi homogeneous curve $x^2 + y^{k+1}$, for some
$k$.  Then we can suppose that $$ f = f_n + f_{n+1} + \cdots =
\sum_{k \geq n}{h_k} = \sum_{k \geq n}{\sum_{i+j=k}{a_{ij} x^i
y^j}}, $$  where $n \geq 3$ and $f_n \neq 0$.

We will write $$ \delta_1 = b_1 \dx + c_1 \dy = \delta_0^{1} +
\delta_1^{1} + \cdots = \sum_{k \geq 0}{\sum_{i+j = k +
1}{(\beta_{ij}^{1} x^i y^j \dx + \gamma_{ij}^{1} x^i y^j \dy)}},
$$ \noindent where the linear part $\delta_0^{1}$ is $(x y) A_0
(\dx \dy)^t$, and $A_0$ is a matrix $2 \times 2$ with complex
coefficients.

If $A_0 = 0$, we have finished. Otherwise, the possibilities of
the Jordan form of $A_0$ are $$ A_0 = \left(
\begin{array}{cc}
\lambda_1 & 0 \\ 0 & \lambda_2
\end{array}
\right), \ \ \ A_0 = \left(
\begin{array}{cc}
\lambda_1 & 0 \\ 1 & \lambda_1
\end{array}
\right). $$ \noindent As $\delta_1$ is not an Euler vector
(because $f$ is not quasi homogeneous), we deduce:
\begin{itemize}
\item If we take the first Jordan form, then (see
the cited demonstration of \cite{4T}) $f_n = x^p y^q$ y $\delta_0
= q x \dx - p y \dy$.  After a sequence of changes of coordinates
we have that $f = x^p y^q$ with $p + q = n \geq 3$, that
contradicts that $f$ is reduced.
\item For the second Jordan form with $\lambda_1 \neq 0$, it has to be
$f_n = 0$, that contradicts that $f$ has its initial part of grade
$n$.
\item For the second option with $\lambda_1 = 0$ we have
$\delta_0^{1} = y \dx$ and, in this situation, the linear of $b_1$
is $y$.  If we precisely apply $\dx$, we obtain 0.
\end{itemize}

In a similar way, you prove the same for $a_2 +  \alpha_1$.
\end{proof}

\vv

\begin{thm}{\label{nocero}}
The natural morphism $\widetilde{M}^{log} \stackrel{\psi}\lra
{\mathcal O}[\frac{1}{f}]$ is an isomorphism if and only if $f$ is
a quasi homogeneous (reduced) curve.
\end{thm}
\begin{proof}
As we pointed, if $f$ is quasi homogeneous then
$\widetilde{I}^{log} = Ann_{\DD}(1/f)$ and therefore $\psi$ is an
isomorphism. Reciprocally, if $\psi$ is an isomorphism, then
$\Ext_{\DD}^2 ({\mathcal O}[1/f], {\mathcal O}) \simeq
\Ext_{\DD}^2 (\widetilde{M}^{log}, {\mathcal O})$. Because of a
result of \cite{MK2}, we have $\Ext_{\DD}^2 ({\mathcal O}[1/f],
{\mathcal O}) = 0$ and, if we take into account proposition
\ref{ext_no_cero}, we obtain that $f$ has to be quasi homogeneous.
\end{proof}

\vv

\noindent {\it Remark.-} The following result can be obtain using
\cite{Torrelli}: if $f$ is not quasi-homogeneous curve then
$Ann_{\mathcal D} (1/f)$ could not be generated by elements of
degree one in $\partial$ and then $Ann_{\mathcal D} (1/f) \neq
\widetilde{I}^{log}$.

\vv

Let us give a new ``differential" proof of a version of the {\em
Logarithmic Comparison Theorem} \cite{4T}.

\begin{thm}{\label{comparlogar}}
The complexes $\Omega^{\bullet} (log f)$ and
$\Omega^{\bullet}[1/f]$ are isomorphic in the correspondent
derived category if and only if $f$ is quasi homogeneous.
\end{thm}
\begin{proof}
If $f$ is quasi homogeneous we have pointed yet that
$\widetilde{M}^{log}$ is isomorphic to ${\mathcal
O}[\frac{1}{f}]$. By the proposition \ref{dual} $(M^{log})^\star
\simeq \widetilde{M}^{log}$ and then we have $$ \Omega^{\bullet}
(log f) \simeq Sol (M^{log}) \simeq DR ((M^{log})^\star) \simeq
DR(\widetilde{M}^{log}) \simeq \Omega^{\bullet} [1/f], $$ where
the first isomorphism is obtained in \cite{Cald} (see also
\cite{Cald3}) and the second one could be found in \cite{MK2}.
Reciprocally, if $f$ is not quasi homogeneous then
${\widetilde{M}}^{log} \not\!{\simeq} \ {\mathcal O}[1/f]$ and, as
both are regular holonomic, neither their De Rham complexes are
isomorphic, that is $$ DR({\widetilde{M}}^{log}) \not\!{\simeq} \
\Omega^{\bullet}[1/f], $$ \noindent using the Riemann-Hilbert
correspondence of Mebkhout-Kashiwara.
\end{proof}

\section{Example in a constructive way of logarithmic comparison in
surfaces.}

\noindent We illustrate in this section an interesting example in
dimension 3 of the situation presented in theorem
\ref{comparlogar} for curves. We consider the surface (see
\cite{Cald}) $h = 0$ with $$ h = x y (x + y)(x z + y), $$
\noindent which is not locally quasi-homogeneous. We prove that:
\begin{itemize}
\item $Ann_{\DD}(1/h) = \widetilde{I}^{log}$.
\item $\widetilde{M}^{log} \simeq (M^{log})^\star$
\end{itemize}
\noindent and we conclude the logarithmic comparison theorem holds
in this case. Although this example appears in \cite{4T}, here the
treatment is under an effective point of view.

We can compute a basis of Der $(\log h)$ with a set of generators
of the syzygies among $h, \frac{\partial h}{\partial_x},
\frac{\partial h}{\partial_y}, \frac{\partial h}{\partial_z}$.  We
obtain
\begin{itemize}
\item $\delta_1 = x \dx + y \dy$,
\item $\delta_2 = x z \dz + y \dz$,
\item $\delta_3 = x^2 \dx - y^2 \dy - x z \dz - y z \dz$,
\end{itemize}
\noindent with $$ \delta_1 (h) = 4 h, \ \ \ \delta_2 (h) = x h, \
\ \ \delta_3 (h) = (2 x - 3 y) h, $$ \noindent and $$ \left|
\begin{array}{ccc}
x & y & 0 \\ 0 & 0 & x z + y \\ x^2 & - y^2 & - x z - y z
\end{array}
\right| = h. $$ \noindent As a multiple of the $b$-function of $h$
in $\DD$ is $$ b(s) = (4 s + 5)(2s + 1) (4 s + 3)(s + 1)^3, $$ and
this polynomial has no integer roots smaller than $-1$, we can
assure that $$ {\mathcal O} [* h] \simeq {\mathcal D}\frac{1}{h}\
. $$ It is easy to check that $Ann_\DD (1/h)$ is equal to
$\widetilde{I}^{log}$. The computations of the $b$-function and
the annihilating ideal of $h^s$ have been made using the
algorithms of \cite{O}, implemented in \cite{Kan}.

\vv

\noindent We calculate (using Gr\"obner bases) a free resolution
of the module ${\mathcal D}/I^{\log}$ where $I^{\log} = (\delta_1,
\delta_2, \delta_3)$ (see \cite{C1}). The first module of syzygies
is generated (in this case) by the relations deduced from the
expressions of the $[\delta_i, \delta_j]$ with $i\neq j$ :

\begin{itemize}
\item $[\delta_1, \delta_2] = \delta_2$,
\item $[\delta_1, \delta_3] = \delta_3$,
\item $[\delta_2, \delta_3] = - x \delta_2$.
\end{itemize}

\noindent The second module of syzygies is generated by only one
element ${\bf s} = (s_1, s_2, s_3)$:
\begin{itemize}
\item $s_1 = - y^2 \dy + x^2 \dx - z y \dz - z x \dz - x$,
\item $s_2 = - y \dz - x z \dz$,
\item $s_3 = y \dy + x \dx - 2$.
\end{itemize}
\noindent The above calculations provide a free resolution of
$M^{log}$.  With a procedure similar to the used in \ref{dual} we
obtain that $(M^{log})^\star$ is the left $\DD$-module associated
to the right $\DD$-module $\DD/(s_1, s_2, s_3)\DD$. Then $$
(M^{log})^\star \simeq \DD/(s_1^t, s_2^t, s_3^t). $$ \noindent It
is enough to compute $s_1^t, s_2^t, s_3^t$ and check (using
Gr\"obner basis) that they span $\widetilde{I}^{\log}$. Hence $$
(M^{log})^\star = (\DD/{\rm Der(log}(h))^\star \simeq \DD/
\widetilde{I}^{log} = \widetilde{M}^{log}. $$

At this point we have obtained that $$ Sol (M^{log}) \simeq DR
((M^{log})^\star) \simeq DR(\widetilde{M}^{log}) \simeq
\Omega^{\bullet} [1/h]$$ where the last two isomorphism are due to
our computations (the first was used in the proof of
\ref{comparlogar}). Taking into account that $\Omega^\bullet($log$
h )\simeq Sol (M^{log})$ (it was showed for this example in
\cite{Cald},\cite{Cald3}) we can deduce that the logarithmic
comparison theorem (i.e. $\Omega^\bullet($log$ h) \simeq
\Omega^\bullet[1/h]$,  \cite{4T})  holds without the ``locally
quasi homogeneous" hypothesis. It is interesting to remark too
(see \cite{Cald}) that $\{ \sigma^F (\delta_1),\sigma^F
(\delta_2), \sigma^F (\delta_3)\}$ do not form a regular sequence
in $\gr^F(\DD)$.  We have $z \eta \zeta - \xi \zeta \notin \langle
\sigma^F (\delta_1), \sigma^F (\delta_2)\rangle$ and $$ (z \eta
\zeta - \xi \zeta) \sigma^F (\delta_3) \in \langle \sigma^F
(\delta_1), \sigma^F (\delta_2)\rangle. $$

\end{document}